\documentclass{amsart}

\usepackage{epsfig}

\usepackage{amssymb}
\usepackage{amsmath}

\newtheorem{thm}{Theorem}[subsection]  
\newtheorem{lem}[thm]{Lemma}	       
      
\newtheorem{prp}[thm]{Proposition}     

\newtheorem{rem}[thm]{Remark}

\makeatletter
\def\theequation{\@arabic\c@equation}
\def\thethm{\@arabic\c@thm}
\def\thelem{\@arabic\c@thm}
\def\thecrlr{\@arabic\c@thm}
\def\theprp{\@arabic\c@thm}
\def\therem{\@arabic\c@thm}
\makeatother

\def\R{{\mathbb{R}}}
\def\NN{{\mathbb{N}}}
\def\sm{\setminus}

\def\ds{\displaystyle}
\def\bs{\boldsymbol}
\def\ol{\overline}

\def\l{\bs{\ell}}
\def\Ll{\bs{L}}
\def\n{\bs{\nu}}
\def\t{\bs{\tau}}
\def\ps{\bs{\psi}}

\def\tl{\widetilde}

\def\B{{\mathcal{B}}}

\def\D{{\mathcal{D}}}
\def\E{{\mathcal{E}}}
\def\F{{\mathcal{F}}}

\def\H{{\mathcal{H}}}

\def\L{{\mathcal{L}}}
\def\N{{\mathcal{N}}}
\def\P{{\mathcal{P}}}

\def\T{{\mathcal{T}}}

\def\FF{{\mathfrak{F}}}

\begin{document}

\title[${W^{2,p}}$-{\it A~Priori\/} Estimates for the Neutral 
Poincar\'e Problem]{$\bs{W^{2,p}}$-A~Priori Estimates  for the 
Neutral Poincar\'e Problem}

\author[D.K. Palagachev]{Dian K. Palagachev}
  
\address{Dipartimento di Matematica, Politecnico di Bari\\ 
Via E. Orabona, 4, 70125 Bari, Italy}
\email{dian@dm.uniba.it, palaga@poliba.it}

\dedicatory{To the memory of Filippo Chiarenza}



\begin{abstract}
A degenerate oblique derivative problem is studied
for uniformly elliptic operators with low regular coefficients
in the framework of Sobolev's classes $W^{2,p}(\Omega)$ for {\em arbitrary\/}
$p>1.$ The boundary operator is prescribed in terms of a directional derivative with respect to the 
vector field $\l$ that becomes tangential to $\partial \Omega$ at the points of some
non-empty subset $\E\subset \partial \Omega$ and is directed outwards
$\Omega$ on $\partial\Omega\setminus\E.$ Under quite general assumptions of the behaviour of $\l,$
we derive {\it a~priori\/} estimates for the $W^{2,p}(\Omega)$-strong solutions   for any $p\in(1,\infty).$  
\end{abstract}

\maketitle

\section*{Introduction}

The lecture deals with regularity in Sobolev's spaces $W^{2,p}(\Omega),$
$\forall\ p\in(1,\infty),$ of the strong solutions to the oblique derivative problem
\begin{equation}\label{P}
\begin{cases}
{\L}u:=
a^{ij}(x)D_{ij}u=f(x)&\text{a.e.}\ \Omega,\\
{\B}u:=
\partial u/\partial \l =\varphi(x)&\text{on} \ \partial \Omega
\end{cases}
\end{equation}
where $\L$ is a uniformly elliptic operator with low regular coefficients
and  $\B$ is prescribed in terms of a directional derivative with respect to
the
unit vector field $\l(x)=(\ell^1(x),\ldots,\ell^n(x))$ defined
on $\partial \Omega,$ $n\geq3.$ Precisely, we are interested in the Poincar\'e problem
\eqref{P} (cf. \cite{Poi,PoP,Pn}), that is, a situation when $\l(x)$ becomes {\it tangential\/} to
$\partial\Omega$ at the points of a non-empty subset $\E$ of $\partial\Omega.$

From a mathematical point of view, \eqref{P} is {\it not\/} an elliptic
boundary value problem. In fact, it follows from the general PDEs theory that \eqref{P} is a
{\it regular (elliptic)\/} problem {\it if and only if\/}  the
Shapiro--Lopatinskij complementary condition is satisfied which means $\l$ must be
transversal to $\partial\Omega$ when $n\geq3$ and $|\l|\neq0$ as $n=2.$ If
$\l$ is {\it tangent\/} to $\partial\Omega$ then \eqref{P} is a
{\it degenerate\/} problem and new effects occur in contrast to the regular case.
It
turns out that the qualitative properties of \eqref{P} depend on the behaviour of
$\l$ near the set of tangency $\E$  and especially on the way 
the normal component $\gamma\n$ of $\l$ (with respect to the 
outward normal $\n$ to $\partial\Omega$) changes or no its sign on the trajectories of $\l$ 
when these cross $\E.$ The main results
were obtained by H\"ormander \cite{H},
Egorov and Kondrat'ev \cite{EK}, Maz'ya
\cite{Mz}, Maz'ya and Paneah \cite{MzPh}, Melin and Sj\"ostrand \cite{MSj},
Paneah \cite{Pn1}
and
good surveys and details can be found in Popivanov and
Palagachev \cite{PoP} and Paneah \cite{Pn}. The problem \eqref{P} has been
studied
in the framework of Sobolev spaces $H^s(\equiv H^{s,2})$ assuming
$C^\infty$-smooth data
and this naturally involved techniques from the pseudo-differential calculus.

The simplest case arises when $\gamma:=\l\cdot\n,$ even if zero on $\E,$ conserves the
sign
on $\partial\Omega.$ Then $\E$ and $\l$ are of {\it neutral\/} type
(a
terminology coming from the physical interpretation of \eqref{P} in
the theory of Brownian
motion, see \cite{PoP}) and \eqref{P} is a problem of Fredholm type (cf. \cite{EK}).
Assume now that  $\gamma$ changes the sign from ``$-$'' to ``$+$'' in
positive direction along the $\l$-integral curves passing through the points of $\E.$
Then $\l$ is of
{\it emergent\/} type and $\E$ is called {\it attracting\/} manifold. The new
effect appearing now is that the {\it kernel\/} of \eqref{P} is {\it
infinite-dimensional\/} (\cite{H}) and to get a well-posed problem one has
to modify \eqref{P} by prescribing the values of $u$ on $\E$ (cf. \cite{EK}).
Finally, suppose the sign of $\gamma$ changes from ``$+$'' to ``$-$'' along the
$\l$-trajectories. Now $\l$ is of {\it submergent\/} type and $\E$ corresponds to a
{\it repellent\/} manifold. The problem \eqref{P} has {\it
infinite-dimensional cokernel\/} (\cite{H}) and Maz'ya and Paneah
\cite{MzPh}
were the first to propose a relevant modification of \eqref{P} by violating
the
boundary condition at the points of $\E.$ As consequence, a Fredholm problem
arises, but the restriction $u|_{\partial\Omega}$ has a finite jump at
$\E.$ What is the common feature of the degenerate problems, independently of the
type of $\l,$ is that the solution ``loses regularity'' near the set of
tangency from the data of \eqref{P} in contrast to the non-degenerate case when any
solution gains two derivatives from $f$ and one derivative from
$\varphi.$ Roughly speaking, that loss of smoothness depends on the {\it order of
contact\/} between $\l$ and $\partial\Omega$ and is given by
the {\it subelliptic\/} estimates obtained for the solutions of degenerate problems (cf.
\cite{GuS1,GuS2,H,MzPh}). Precisely, if $\l$ has a contact of order $k$ with $\partial\Omega$ then
the solution of \eqref{P} gains $2-{k}/{(k+1)}$ derivatives from $f$ and 
$1-{k}/{(k+1)}$ derivatives from $\varphi.$
\begin{center}
\begin{minipage}[t]{61mm} 
\centerline{\includegraphics[scale=.8]{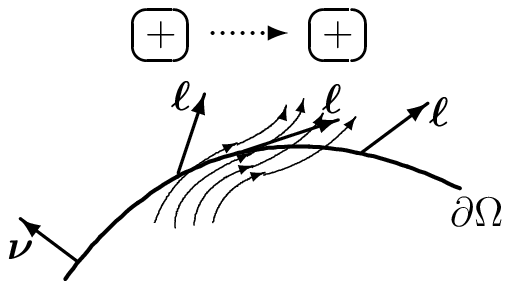}}
\vspace*{0.1cm}
\centerline{{(a)}~neutral vector field $\l$}
\end{minipage}\hfill\begin{minipage}[t]{61mm}
\centerline{\includegraphics[scale=.8]{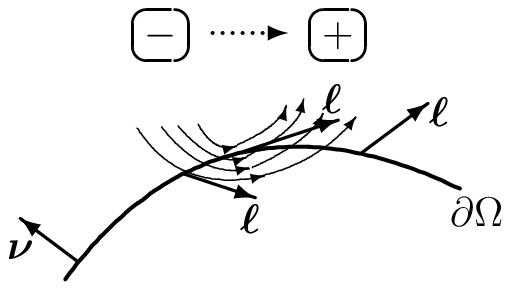}}
\vspace*{0.1cm}
\centerline{{(b)}~emergent vector field $\l$}
\end{minipage}

\vspace*{0.2cm}

\begin{minipage}{61mm}
\centerline{\includegraphics[scale=.8]{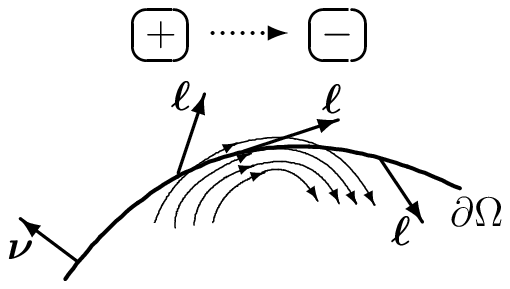}}
\vspace*{0.1cm}
\centerline{{(c)}~submergent vector field $\l$}
\end{minipage}
\end{center}
For what concerns the geometric structure of $\E,$ it was supposed initially to be a
submanifold of $\partial\Omega$ of codimension one. Melin and Sj\"ostrand
\cite{MSj} and Paneah \cite{Pn1} were the first to study the Poincar\'e problem \eqref{P}
in a
more general situation when $\E$ is a massive subset of $\partial\Omega$ with positive surface measure,
 allowing $\E$ to contain arcs of $\l$-trajectories of {\it finite\/} length. Their results
were extended by Winzell (\cite{W1,W2}) to the framework of H\"older's spaces who studied  \eqref{P}
assuming $C^{1,\alpha}$-smoothness of the coefficients
of $\L.$ It is worth noting that $\l$ has automatically an {\it infinite\/}
order of contact with $\partial\Omega$ when $\E$ is a massive subset of the boundary.

To deal with non-linear Poincar\'e problems, however,
we have to dispose of precise information on the linear problem \eqref{P} with
coefficients less regular than $C^\infty$ (see \cite{Pl0,PoK0,PoK,PoP}). Indeed, {\it a~priori\/} 
estimates in $W^{2,p}$ for solutions to \eqref{P} would imply easily  pointwise
estimates for $u$ and $Du$ for suitable values of $p>1$ through the Sobolev imbeddings.
This way, we are naturally led to consider the problem \eqref{P}
in a {\it strong\/} sense, that is, to searching for solutions lying in
$W^{2,p}$ which satisfy  $\L u=f$ almost everywhere (a.e.) in $\Omega$
and $\B u=\varphi$ holds in the sense of trace on $\partial\Omega.$

In the papers \cite{GuS1,GuS2} by Guan and Sawyer
solvability and precise subelliptic estimates have been obtained for
\eqref{P} in $H^{s,p}$-spaces ($\equiv W^{s,p}$
for integer $s!$). However, \cite{GuS1} treats operators with
$C^\infty$-coefficients and this determines the technique involved and the results obtained, 
while in \cite{GuS2} the coefficients are
$C^{0,\alpha}$-smooth, but the field $\l$ is of finite type, that is, it has a
{\it finite\/} order of contact with $\partial\Omega.$

The main goal of this lecture is to derive {\it a~priori\/} estimates 
in Sobolev's classes $W^{2,p}(\Omega)$ with {\it any\/} $p\in(1,\infty)$
for the solutions of the Poincar\'e problem \eqref{P},  
weakening both Winzell's assumptions on $C^{1,\alpha}$-regularity of the coefficients of $\L$
and these of Guan and Sawyer on the {\it finite type\/} of $\l.$ 
We are dealing  with the simpler case when
$\gamma$ preserves its sign on $\partial\Omega$ which means the field
$\l$ is of {\it neutral type.\/} 
Of course, the loss of smoothness 
mentioned, imposes some more regularity of the data near the set $\E.$
We assume the coefficients of $\L$ to be Lipschitz continuous near $\E$ while only continuity (and even
discontinuity controlled in $VMO$) is allowed away from $\E.$ 
Similarly, $\l$ is a Lipschitz vector field on
$\partial\Omega$ with Lipschitz continuous first derivatives near $\E,$ and {\it no
restrictions\/} on the order of contact with $\partial\Omega$ are required.
Regarding the tangency set 
$\E,$ it may have positive surface measure and is restricted only to a sort of  {\it
non-trapping\/} condition that all trajectories of $\l$ through the points of
$\E$ are non-closed  and leave $\E$ in a finite time. 

The technique adopted is based on a dynamical system approach employing
the fact that $\partial u/\partial\l$ is a local strong solution, near $\E,$ 
to a Dirichlet-type problem with right-hand side depending on the solution $u$ itself. 
Application of the $L^p$-estimates for such problems leads to the functional inequality
\eqref{eqzeta} for suitable $W^{2,p}$-norms of $u$ on a family of subdomains which, starting away from $\E,$
evolve along the $\l$-trajectories and exhaust a sort of their tubular neighbourhoods. Fortunately, that is an inequality
with advanced argument and the desired $W^{2,p}$-estimate follows by iteration with respect to the curvilinear
parameter on the 
trajectories of $\l.$ Another advantage of this approach is the {\it improving-of-integrability\/} property obtained for the
solutions of \eqref{P}. Roughly speaking, it asserts that the problem  \eqref{P}, even if a {\it degenerate\/} one, 
behaves as an {\it elliptic\/} problem for what concerns the degree $p$ of
integrability. In other words, the second derivatives of any solution to \eqref{P} will
have the same rate of integrability as $f$ and $\varphi.$ We refer the reader to the paper \cite{P-IMRN} for outgrowths of 
the $W^{2,p}$-{\it a~priori\/} estimates, such as uniqueness in $W^{2,p}(\Omega),$ $\forall\ p>1,$ of 
the strong solutions to \eqref{P} as well as its Fredholmness.

Concluding this introduction, we should mention the articles \cite{MPV1,MPV2,Pl1} where similar results
 have been obtained by different technique in the 
particular case when the tangency set $\E$ contains trajectories of $\l$ with positive, but {\it small enough\/} lengths.

\section{Hypotheses and the Main Result}\label{s1}

Hereafter $\Omega\subset {\R}^n,$ $n\geq3,$ will be a bounded domain with
reasonably smooth boundary and $\n(x)=\big(\nu^1(x),\,\ldots,\,\nu^n(x)\big)$ stands for the unit
{\it outward\/} normal to $\partial \Omega$ at $x\in\partial \Omega.$ Consider a unit vector field
$\l(x)=\big(\ell^1(x),\ldots,\ell^n(x)\big)$ on
$\partial \Omega$ and let
$\l(x)=\t(x)+\gamma(x)\n(x),$ 
where $\t\colon\partial\Omega\to\R^n$
is the projection of $\l(x)$ on the hyperplane tangent to $\partial \Omega$
at $x\in\partial \Omega$ and $\gamma\colon \partial\Omega\to\R$ is the inner product
$\gamma(x):=\l(x)\cdot\n(x).$  The set of zeroes of $\gamma,$ 
$$
\E:=\big\{x\in\partial \Omega\colon\ \gamma(x)=0\big\},
$$
is indeed the subset of $\partial \Omega$ where the field
$\l(x)$ becomes tangent to it.

Fix $\N\subset\ol\Omega$ to be a closed neighbourhood of $\E$ in
$\ol\Omega.$ We suppose $\L$ is
a uniformly elliptic operator with measurable coefficients, satisfying
\begin{equation}\label{eq1.1}
\lambda^{-1}|\xi|^2\leq a^{ij}(x)\xi_i\xi_j\leq \lambda|\xi|^2\quad
\text{a.a.}\ x\in\Omega,\ \forall \xi\in {\R}^n,\qquad a^{ij}(x)=a^{ji}(x)
\end{equation}
for some positive constant $\lambda.$ Regarding the regularity of the data, we assume
\begin{equation}\label{eq1.2}
\begin{cases}
a^{ij}\in VMO(\Omega)\cap C^{0,1}(\N),\\
\partial\Omega\in C^{1,1},\quad \partial\Omega\cap \N\in C^{2,1},\quad
\ell^i\in C^{0,1}(\partial \Omega)\cap
C^{1,1}(\partial \Omega\cap \N)
\end{cases}
\end{equation}
with $VMO(\Omega)$ being the Sarason class of functions of vanishing mean oscillation and
$C^{k,1}$ denotes the space of functions with Lipschitz continuous $k$-th order
derivatives. Let us point out that \eqref{eq1.1},
 \eqref{eq1.2} and the Rademacher theorem give  $a^{ij}\in L^\infty(\Omega)\cap W^{1,\infty}(\N).$ 
For what concerns the boundary operator $\B,$ we assume
\begin{equation}\label{eq1.3}
\begin{cases}
\gamma(x)=\l(x)\!\cdot\!\n(x)\geq 0\quad  \forall \, x\in\partial \Omega,\quad\ \text{and}\\
\text{the arcs of the\ } \l\text{-trajectories lying in\ } \E\ \text{(which coincide with these of}\ \t)\\ 
\text{are all\em\  non-closed\em\ and of\em\ finite lengths.}
\end{cases}
\end{equation}
The first assumption simply means that $\l(x)$ is either tangential to
$\partial \Omega$ or is directed outwards $\Omega,$ 
that is, the field $\l$ is of {\it neutral
type\/} on $\partial \Omega,$ while the second one is a sort of {\it non-trapping\/} 
condition on the tangency set $\E.$ It implies that
  the $\l$-integral curves {\it leave\/}  $\E$
in a {\it finite time\/} in both directions.
\begin{figure}[hbt] 
\includegraphics[scale=.8]{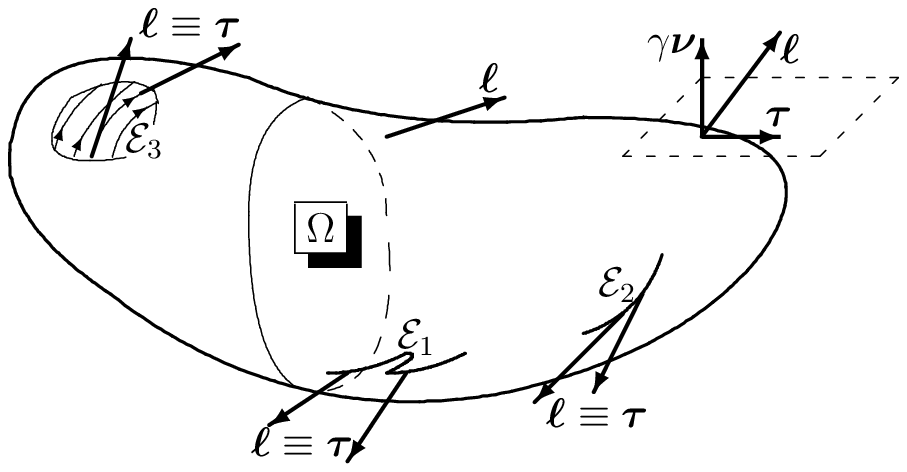}
\caption[]{The set of tangency $\E$ is the union $\E_1\cup\E_2\cup\E_3$ where
$\mathrm{codim\,}_{\partial\Omega}\E_1=\mathrm{codim\,}_{\partial\Omega}\E_2=1$ while
$\mathrm{meas\,}_{\partial\Omega}\E_3>0.$ The vector field $\l$ is transversal to $\E_1$ and tangent to $\E_2.$
Actually, $\E_2$ consists of an arc of $\t$-trajectory, whereas $\E_3$ is union of such arcs.}
\end{figure}

Throughout the text $W^{k,p}$ stands for the Sobolev class of functions with
$L^p$-summable weak derivatives up to order $k\in\NN$ while $W^{s,p}(\partial\Omega)$ with $s>0$
non-integer and $p\in(1,+\infty),$  is the
Sobolev space of fractional order on $\partial\Omega.$ Further, we use the standard
parameterization $t\mapsto \ps_{\Ll}(t;x)$ for the {\it trajectory\/} (equivalently, {\it phase curve,
maximal integral curve\/}) of a
given vector field $\Ll$ passing through a point $x,$ that is, $\partial_t \ps_{\Ll}(t;x)=
\Ll\big(\ps_{\Ll}(t;x)\big)$ and $\ps_{\Ll}(0;x)=x.$

We will employ below an extension of the field $\l$ near $\partial\Omega$ which preserves therein its regularity and 
geometric properties. All the results and proofs in the sequel work for such an {\em arbitrary\/} $\l$-extension but,
in order to make more evident some 
\begin{figure}[htb] 
\includegraphics[scale=.8]{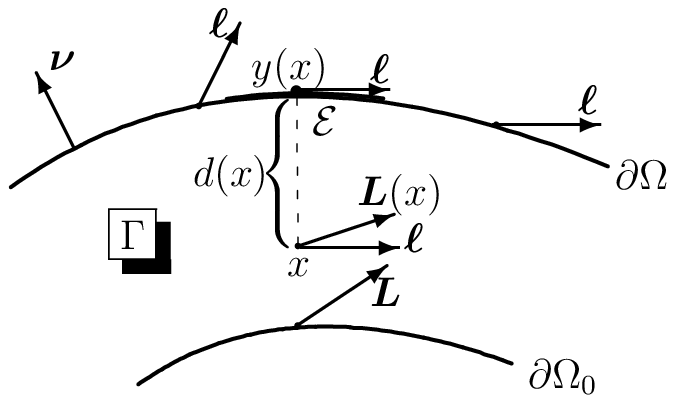}
\caption[]{}
\end{figure}
geometric constructions, we prefer to introduce a {\em special\/} extension as follows. For each $x\in\R^n$
near $\partial\Omega$ set
$d(x)=\mathrm{dist\,}(x,\partial\Omega)$ and define $\Gamma:=\{x\in\R^n\colon\
d(x)\leq d_0\}$ with small $d_0>0.$ Letting
$\Omega_0:=\Omega\sm\Gamma$ and $y(x)\in\partial\Omega$ for the unique point closest to $x\in\Gamma,$ we have
(see~\cite[Chapter~14]{GT})
$y(x)\in C^{0,1}(\Gamma)$ while $y(x)\in C^{1,1}$ near $\E.$ Regarding
the distance function $d(x)=|x-y(x)|,$ it is Lipschitz continuous in $\Gamma$ and inherits
the regularity of $\partial\Omega$ at $y(x)$ when considered on the
parts of $\Gamma$ lying in/out $\Omega,$ but its normal derivative has a finite jump on $\partial\Omega.$ 
Anyway, it is a routine to check $\big(d(x)\big)^2\in C^{1,1}(\Gamma).$ Setting
$\Ll(x)$ for the normalized representative of $\l(y(x))+\big(d(x)\big)^2\n(y(x))$ $\forall x\in \Gamma,$
it results $|\Ll(x)|=1,$ $\Ll|_{\partial\Omega}=\l,$ $\Ll|_{\E}=\t$ and $\Ll\in C^{0,1}(\Gamma)\cap
C^{1,1}(\Gamma\cap \N).$ Moreover, the field $\Ll$ is strictly transversal to $\partial\Omega_0.$

As consequence of the non-trapping condition \eqref{eq1.3}, the compactness of $\E$ and the
semi-continuity properties of the lengths of the $\t$-maximal integral curves, it is not hard to get that (see~\cite[Proposition~3.1]{W2} and
\cite[Proposition~3.2.4]{PoP}) {\it under the hypotheses
\eqref{eq1.2} and \eqref{eq1.3}, there is a finite upper bound
$\kappa_0$ for the arclengths of the $\t$-trajectories lying in $\E.$ Moreover, each point of
$\Gamma$ can be reached from $\partial\Omega_0$ by an $\Ll$-integral curve of length at most
$\kappa=\mathrm{const}>0.$}

In what follows, the letter $C$ will denote a generic
constant depending on known quantities defined by the data of \eqref{P}, that is, on $n,$ $p,$ $\lambda,$
the respective norms of the coefficients of $\L$ and $\B$ in $\Omega$ and $\N,$  the regularity of $\partial\Omega$
and the constants $\kappa_0$  and $\kappa.$

In order to control precisely the regularity of $u$ near the tangency set $\E,$ we have to introduce the
appropriate functional spaces. For, take an arbitrary $p\in(1,\infty)$ and define the Banach spaces
$$
\F^p(\Omega,\N):=\left\{f\in L^p(\Omega)\colon\ \partial f/\partial \Ll\in
L^p(\N) \right\}
$$
equipped with norm $\|f\|_{\F^p(\Omega,\N)}:=
\|f\|_{L^p(\Omega)}+\|\partial f/\partial \Ll\|_{L^p(\N)},$
and 
$$
\Phi^p(\partial\Omega,\N)
:=\left\{\varphi\in W^{1-1/p,p}(\partial\Omega)\colon\ \varphi\in
		   W^{2-1/p,p}(\partial\Omega\cap\N)\right\}
$$
normed by $\|\varphi\|_{\Phi^p(\partial\Omega,\N)}:=
\|\varphi\|_{W^{1-1/p,p}(\partial\Omega)}+
\|\varphi\|_{W^{2-1/p,p}(\partial\Omega\cap\N)}.$

Our main result asserts that the couple $(\L,\B)$ {\it improves the integrability\/} of solutions to \eqref{P} for any $p$ in the range 
$(1,\infty)$ and, moreover, provides for an {\it a~priori estimate\/} in the
$L^{p}$-Sobolev scales for any such solution.
\begin{thm}\label{thmINC}
Under the hypotheses \eqref{eq1.1}--\eqref{eq1.3} let
$u\in W^{2,p}(\Omega)$ be a strong solution of the problem \eqref{P}
with $f\in \F^q(\Omega,\N)$ and
$\varphi\in \Phi^q(\partial\Omega,\N)$ where $1<p\leq q<\infty.$

Then $u\in W^{2,q}(\Omega)$  and there is an absolute constant $C$ such that
\begin{equation}\label{eq2.4}
\|u\|_{W^{2,q}(\Omega)}  \leq C\left(\|u\|_{L^{q}(\Omega)}
+ \|f\|_{\F^q(\Omega,\N)} +
\|\varphi\|_{\Phi^q(\partial\Omega,\N)}\right).
\end{equation}
\end{thm}
Let us point out reader's attention that the directional derivative $\partial u/\partial\Ll$
of each $W^{2,p}$-solution to \eqref{P} belongs to $W^{2,p}(\N).$ For, $\partial u/\partial\Ll\in W^{1,p}(\N)$ 
and taking the difference quotients in \eqref{P} in the direction of $\Ll$ 
(cf. \cite[Chapter~8 and Lemma~7.24]{GT}) gives that $\partial u/\partial\Ll\in W^{2,p}(\N)$ is a strong local solution to the Dirichlet problem
\begin{equation}\label{eqderDP}
\begin{cases}
\L \left(\frac{\partial u}{\partial \Ll}\right)=\frac{\partial f}{\partial \Ll}+2a^{ij}D_jL^kD_{ki}u+a^{ij}D_{ij}L^kD_ku
 - \frac{\partial a^{ij}}{\partial \Ll}D_{ij}u
\qquad\text{a.e.}\ \N,\\
\frac{\partial u}{\partial \Ll}=\varphi \qquad\text{on}\ \partial\Omega\cap\N
\end{cases}
\end{equation}
where $\Ll(x)=(L^1(x),\ldots,L^n(x))\in C^{1,1}(\N).$ Therefore, once having proved $u\in W^{2,q}(\Omega)$ and the estimate 
\eqref{eq2.4}, we have
$$
\|\partial u/\partial\Ll\|_{W^{2,q}(\tl\N)}  \leq C'\left(\|u\|_{L^{q}(\Omega)}
+ \|f\|_{\F^q(\Omega,\N)} +
\|\varphi\|_{\Phi^q(\partial\Omega,\N)}\right)
$$
for any {\em closed\/} neighbourhood $\tl\N$ of $\E$ in $\ol\Omega,$ $\tl\N\subset\N,$ by means of the $L^p$-theory of 
uniformly elliptic equations (see \cite{ADN} or \cite[Chapter~9]{GT}).
In other words, {\it if a strong solution $u$ to \eqref{P} belongs to $W^{2,q}(\Omega)$ then
$\partial u/\partial\Ll\in W^{2,q}(\N)$ automatically, provided $f\in \F^q(\Omega,\N)$ and
$\varphi\in \Phi^q(\partial\Omega,\N).$}

\section{Proof of Theorem \ref{thmINC}}

Fix hereafter $\N'\subset\N''\subset\N$ to be closed neighbourhoods of $\E$ in
$\ol{\Omega}$ with $\N''$ so ``narrow'' that $\N''\subset
\Omega\sm \Omega_0$ (see Figure~3£). The next result is an immediate consequence of $\gamma(x)>0$
$\forall x\in \partial\Omega\sm\N'$ and
 the $L^p$-theory of {\it regular\/} oblique derivative problems for uniformly elliptic operators with $VMO$ principal
coefficients (cf. \cite[Theorem~2.3.1]{MPS}).
\begin{prp}\label{lemA1}
Assume \eqref{eq1.1}, \eqref{eq1.2} and $\gamma(x)>0$ $\forall x\in\Omega\setminus\E,$ and let
$u\in W^{2,p}(\Omega)$ be a solution to \eqref{P}
with $f\in L^q(\Omega)$ and $\varphi\in W^{1-1/q,q}(\partial\Omega),$ where $1<p\leq q<\infty.$

Then $u\in W^{2,q}(\Omega\sm \N')$ and there is a constant such that
\begin{equation}\label{eq2.5}
\|u\|_{W^{2,q}(\Omega\sm\N')} \leq C\left(
\|u\|_{L^{q}(\Omega)}+
\|f\|_{L^{q}(\Omega)}+
\|\varphi\|_{W^{1-1/q,q}(\partial\Omega)}\right).
\end{equation}
\end{prp}

To derive the improving-of-integrability near the tangency set $\E,$
we consider any solution of the problem \eqref{P} for which 
$a^{ij},\ \partial a^{ij}/\partial\Ll\in L^\infty(\N)$ in view of 
\eqref{eq1.2}\footnote{It will be clear from the considerations given below 
that instead of Lipschitz continuity
of the coefficients of $\L$ in $\N$ as \eqref{eq1.2} asks, it suffices to
have essentially bounded their directional derivatives with respect to the field $\Ll.$} and
$f,$ $\partial f/\partial\Ll\in L^q(\N)$ and
$\varphi\in W^{2-1/q,q}(\partial\Omega\cap\N)$ by hypotheses.
\begin{lem}\label{lemA2}
Under the assumptions of Theorem~$\ref{thmINC},$ the solution $u$ of
$\eqref{P}$ belongs to $u\in W^{2,q}(\N'')$ and there is a constant such that
\begin{equation}\label{eq2.6}
\|u\|_{W^{2,q}(\N'')} \leq C\left(
\|u\|_{L^{q}(\Omega)}+
\|f\|_{\F^{q}(\Omega,\N)}+
\|\varphi\|_{\Phi^{q}(\partial\Omega,\N)}\right).
\end{equation}
\end{lem}
\begin{proof}
Take an arbitrary point $x_0\in\E.$ 
According to \eqref{eq1.3}, the $\Ll$-trajectory through $x_0$ leaves
 $\E$ in both directions for a finite time, that is, 
$\ps_{\Ll}(t^-;x_0)\in \N''\sm\N',$
$\ps_{\Ll}(t^+;x_0)\in \R^n\sm\ol\Omega$ (see Figure~3£) for suitable $t^-<0<t^+.$ 
\begin{figure}[htb] 
\includegraphics[scale=.8]{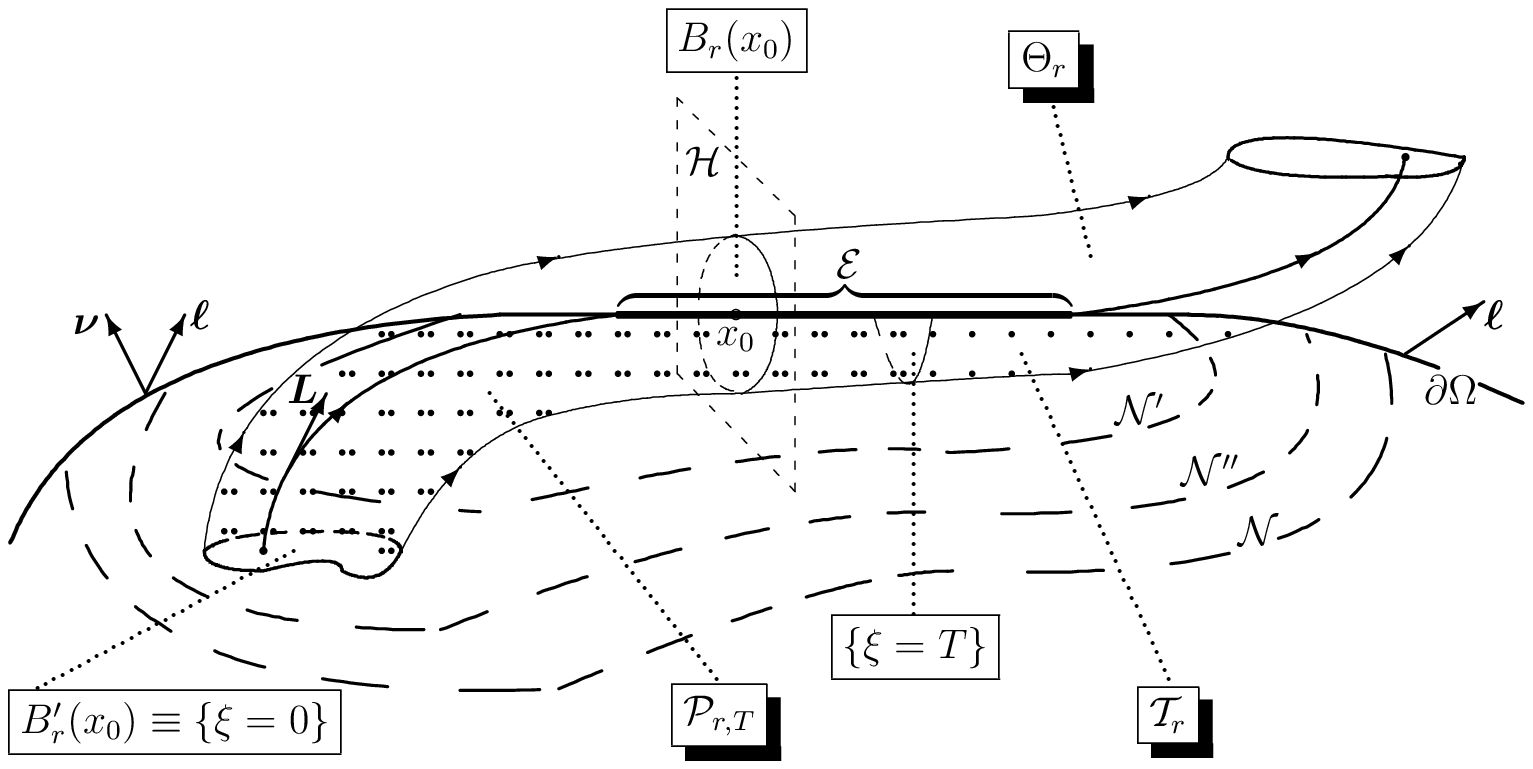}
\caption[]{$\T_r$ is the dotted set, while the double-dotted one is $\P_{r,T}.$}
\end{figure}

Set $\H$ for the $(n-1)$-dimensional
hyperplane through $x_0$ and orthogonal to $\Ll(x_0),$ and define
$$
B_r(x_0):=\left\{x\in\H\colon\quad |x-x_0|<r\right\}
$$
with $r>0$ to be chosen later. It follows from the Picard inequality\footnote{$|\ps_{\Ll}(t;x')-\ps_{\Ll}(t;x'')|\leq \text{\sl
e\,}^{t\|\Ll\|_{C^1(\N)}}|x'-x''|$ for all $x',\ x''\in\N.$} that if $r$ is small enough, then the flow of $B_r(x_0)$ along the $\Ll$-trajectories at time $t^-,$
$$
B'_r(x_0):=\ps_{\Ll}(t^-;B_r(x_0)):=\left\{\ps_{\Ll}(t^-;y)\colon\quad
y\in B_r(x_0)\right\}
$$
is {\it entirely contained\/} in $\N''\sm\N'$ whence
$B'_r(x_0)\cap\E=\emptyset.$ The set
$$
\Theta_r:= \left\{\ps_{\Ll}(t;x')\colon\quad x'\in B'_r(x_0),\quad
t\in(0,t^+-t^-) \right\}
$$
is an $n$-dimensional neighbourhood of the $\Ll$-trajectory
through $x_0$  and defining
$$
\T_r:=\Theta_r\cap\Omega,
$$
the boundary $\partial\T_r$ is composed of
the ``base'' $B'_r(x_0)$ and the ``lateral'' components
$\partial_1\T_r:=\partial\T_r\cap\partial\Omega$ and
$\partial_2\T_r:=(\partial\T_r\cap \Omega)\sm B'_r(x_0).$ Indeed, $\T_r\subset \N''$ if $r>0$ is small enough.

We will derive \eqref{eq2.6} in $\T_r$ after that the desired estimate will follow by covering the
compact $\E\subset\partial\Omega$ by a finite number of sets like $\ol{\T_r}.$  Our strategy is based
on a representation of $u(x)$ in $\T_r$ by means of $u(x')$ with $x'=\ps_{\Ll}(-\xi(x);x)\in B'_r(x_0)$
for some $\xi(x)>0,$ and the integral of $\partial u/\partial\Ll$ along the $\Ll$-trajectory joining
$x'$ with $x.$ Thus the Sobolev norm of $u$ will be expressed by
the respective norm of $\partial u/\partial\Ll$ and that of $u$ itself near
$B'_r(x_0)$ where we dispose of \eqref{eq2.5}. Concerning $\partial
u/\partial\Ll,$  it is a local solution of Dirichlet problem near $\E$ with
right-hand side depending on $u.$

Let $\mu\colon \H\to\R^+$ be a $C^\infty$ cut-off function such that
\begin{equation}\label{eqMU}
\mu(y)=\begin{cases}
1 & y\in B_{r/2}(x_0),\\
0 & y\in \H\sm B_{3r/4}(x_0)
\end{cases}
\end{equation}
and extend it to $\R^n$ as constant on the $\Ll$-trajectory
through $y\in\H.$ The function $U(x):=\mu(x)u(x)$
is a $W^{2,p}(\N)$-solution of 
\begin{equation}\label{eqUeq}
\begin{cases}
{\L}U=F(x):=\mu f+2a^{ij}D_j\mu D_iu +u a^{ij}D_{ij}\mu\qquad
	 \text{a.e.}\ \T_r,\\[2pt]
\partial U/\partial\Ll=\Phi:=
\begin{cases}
\mu\varphi & \text{on}\ \partial_1\T_r,\\
0	   & \text{near}\ \partial_2\T_r,\\
\mu\partial u/\partial\Ll  & \text{on}\ B'_r(x_0)\subset \N''\sm\N'.
\end{cases}
\end{cases}
\end{equation}
Indeed, $u\in W^{2,p}(\N)$ implies $Du\in L^{np/(n-p)}$ if $p<n$ and $Du\in L^s$
$\forall s>1$ when $p\geq n,$ whence $F\in L^{q'}(\N)$ with
\begin{equation}\label{eqq'}
q':=
\begin{cases}\ds
\min\left\{q,\frac{np}{n-p}\right\} & \text{if}\ p<n,\\
q & \text{if}\ p\geq n.
\end{cases}
\end{equation}
Further, $\partial F/\partial\Ll\in L^{q'}(\N'')$ as consequence of
\eqref{eqderDP}, $\partial u/\partial\Ll\in W^{2,q}(\N''\sm\N')$
by Proposition~\ref{lemA1} whence $\Phi\in
W^{2-1/q,q}(\partial\T_r).$ Thus \eqref{eq1.1}, \eqref{eq1.2}, $\T_r\subset\N''$ and \eqref{eqderDP} give that
$$
V(x):=\partial U/\partial\Ll 
$$
is a  $W^{2,p}(\T_r)$-solution of the Dirichlet problem
\begin{equation}\label{eq2.7}
\begin{cases}
{\L}V=\partial F/\partial\Ll+2 a^{ij}D_{j}L^k D_{ik}U+
a^{ij}D_{ij}L^kD_kU -\frac{\partial a^{ij}}{\partial\Ll}D_{ij}U
 \quad\text{a.e.}\ \T_r,\\
V=\Phi \quad\text{on} \ \partial \T_r.
\end{cases}
\end{equation}
Now we pass from $x\in\Theta_r$ into the new variables $(x',\xi)$ with $x'=\ps_{\Ll}(-\xi(x);x)\in B'_r(x_0)$
and $\xi\colon\ \Theta_r\to(0,t^+-t^-),$
$\xi(x)\in C^{1,1}(\Theta_r).$ The transform $x\mapsto (x',\xi)$ defines a
$C^{1,1}$-diffeomorphism because the field $\Ll$ is transversal to
$B'_r(x_0).$ Moreover, 
$\partial/\partial\Ll\equiv \partial/\partial\xi,$
$\ps_{\Ll}(t;x')=(x',t)$ and  $V(x',\xi)=\partial
U(x',\xi)/\partial\xi$ as $(x',\xi)\in\T_r.$ Since $V(x',\xi)$ is an absolutely continuous function in
$\xi$ for a.a. $x'\in B'_r(x_0))$ (after redefining it, if necessary, on a set of zero measure)
we get
\begin{equation}\label{eqIFF}
U(x',\xi)   = U(x',0)+\int_0^\xi V(x',t)dt\qquad \text{for a.a.}\
(x',\xi)\in \T_r,
\end{equation}
where the point $(x',0)\in B'_r(x_0)$  lies in $\N''\sm\N'$ and $U(x',0)\in W^{2,q}$ there by Proposition~\ref{lemA1}, 
the Fubini theorem and \cite[Remark~2.1]{Pl2}. Passing to the new variables $(x',\xi)$ in \eqref{eq2.7}, taking the derivatives
of \eqref{eqIFF} up to second order and substituting them into the right-hand
side of \eqref{eq2.7}, this last  reads
\begin{equation}\label{eq2.8f}
\begin{cases}
{\L}'V=F_1(x',\xi) +\ds\int_0^\xi \D_2(\xi) V(x',t)dt
		      &\text{a.e.}\ \T_r,\\
V=\Phi &\text{on} \ \partial \T_r,
\end{cases}
\end{equation}
where ${\L}'$ is the operator $\L$ in terms of $(x',\xi)=(x'_1,\ldots,x'_{n-1},\xi),$ 
\begin{align}\label{eqF1}
\nonumber
F_1(x',\xi):=&\ \partial F/\partial\Ll + \D_1 V(x',\xi)+\D'_1 U(x',\xi)+\D'_2 U(x',0),\\[-8pt]
 & \\[-8pt]
\nonumber
\D_2(\xi) V(x',t) :=&\ \sum_{i,j=1}^{n-1} A^{ij}(x',\xi)D_{x'_ix'_j}V(x',t),\qquad
A^{ij}\in L^\infty,
\end{align}
$\D_1,$ $\D'_1,$ $\D'_2$ are linear differential operators with
$L^\infty$-coefficients, $\mathrm{ord\,}\D_1=\mathrm{ord\,}\D'_1=1,$ $\mathrm{ord\,}\D'_2=2.$
The Sobolev imbedding theorem implies $F_1\in L^{q'}(\T_r)$ with $q'$ given by
\eqref{eqq'} as consequence of $\partial F/\partial\Ll \in L^{q'}(\N''),$
$U(x',0)\in W^{2,q}(B'_r(x_0))$ and $U,\ V\in W^{2,p}(\N'').$ Nevertheless
 the second-order operator $\D_2(\xi)$ has a quite
rough characteristic form which is neither symmetric nor sign-definite,
the improving-of-integrability holds for \eqref{eq2.8f} thanks to the
particular structure of $\T_r$ as union of $\Ll$-trajectories through
$B'_r(x_0).$ Actually, we will show that if $V\in W^{2,q'}$ on a subset of
$\T_r$ with $\xi<T,$ then $V$ remains a $W^{2,q'}$-function on a larger subset
with $\xi<T+r$ for small enough $r,$ after that the higher integrability of $U$
will follow from Proposition~\ref{lemA1} and \eqref{eqIFF}. For, take an arbitrary
$T\in(0,t^+-t^-)$  and define
$$
\P_{r,T}:=\left\{(x',\xi)\in\T_r\colon\quad \xi<T\right\}.
$$
For a fixed $r>0,$ $\big\{\P_{r,T}\big\}_{T\geq0}$ is a non-decreasing family of
domains exhausting $\T_r$ and $\P_{r,T}\equiv\T_r$ for values of $T$ greater than
the {\it maximal exit-time\/}
$$
T_{\max}:=\sup_{x'\in B'_r(x_0)}\sup\left\{t>0\colon\quad \ps_{\Ll}(t;x')\in\Omega,\ x'\in B'_r(x_0)\right\}.
$$
\begin{prp}\label{prpIMP}
Let $T\in(0,t^+-t^-)$ and consider the solution $V\in W^{2,p}(\T_r)$
of the problem $\eqref{eq2.8f}.$ Suppose $V\in W^{2,q'}(\P_{r,T})$ where $q'$ is given by $\eqref{eqq'}.$

There exists an $r_0>0$ 
such that $V\in W^{2,q'}(\P_{r,T+r})$ for all $r<r_0.$
\end{prp}
\begin{proof}
There are three possible cases to be distinguished.
\paragraph*{\it Case A: $T+3r<T_{\max}$} We have $\P_{r,T}\subset \P_{r,T+3r}\subset\T_r\equiv \P_{r,T_{\max}}$ and
consider the $C^\infty$-function
$\eta\colon\ \R\to[0,1]$  such that
\begin{equation}\label{eqeta}
\eta(\xi)=\begin{cases}
1  & \text{as}\ \xi\in(-\infty,T+r],\\
\text{strictly decreases} & \text{as}\ \xi\in (T+r,T+2r),\\
0 & \text{as}\ \xi\geq T+2r.
\end{cases}
\end{equation}
Setting $\tl V(x',\xi):=\eta(\xi)V(x',\xi),$
it follows $\L'\tl V=\eta(\L'V)+\L_1 V$ where $\L_1$ is a first-order differential operator
with $L^\infty$-coefficients depending on these of $\L'$ and on the derivatives of $\eta.$
Therefore,
\begin{align}\label{eq2.9}
\L'\tl V=&\ \eta F_1+ \L_1V +\eta(\xi) \int_0^\xi \D_2(\xi) V(x',t)dt\\
\nonumber
	=&\ \eta F_1 +\L_1V+\int_0^\xi \frac{\eta(\xi)}{\eta(t)} \D_2(\xi)\tl
V(x',t)dt
\end{align}
because $\D_2(\xi)$ is a second-order operator acting in the $x'$-variables only.

We set $\Omega_r\subset \P_{r,T+3r}\sm \P_{r,T-3r}$ for a $C^{1,1}$-smooth domain containing 
$\P_{3r/4,T+2r}\sm \P_{3r/4,T-2r}$
and such that 
$$
r^{-1}\Omega_r:=\left\{(\tl y',\tl \xi)\colon\quad \tl y'=x'/r,\ \tl \xi=(\xi-T)/r,\
(x',\xi)\in\Omega_r\right\} \in C^{1,1}
$$
uniformly in $r.$ 
The boundary $\partial\Omega_r$ consists of the ``lateral'' parts $\partial_1\Omega_r:=\partial\Omega_r\cap\partial\Omega$
and $\partial_2\Omega_r:=\partial\Omega_r\cap \Omega\cap \{\xi\in(T-2r,T+2r)\}
\subset \big(\P_{r,T+2r}\sm\P_{r,T-2r}\big)\sm \big(\P_{3r/4,T+2r}\sm\P_{3r/4,T-2r}\big),$
and of two $C^{1,1}$-smooth components $\partial\Omega_r^\pm$ lying
in $\P_{r,T+3r}\sm\P_{r,T+2r}$ and  $\P_{r,T-2r}\sm\P_{r,T-3r},$
 respectively. The properties of $\mu$ (cf. \eqref{eqMU}) ensure
$U\equiv0,$ $V\equiv0,$ $\tl V\equiv 0$ on $\T_r\sm\T_{3r/4}$ whence
$\tl V\equiv 0$ near $\partial_2\Omega_r.$ 

For an arbitrary $(x',\xi)\in\Omega_r,$ the factor $\eta(\xi)/\eta(t)$
in \eqref{eq2.9} vanishes when $\xi\geq T+2r$ while
$\eta(\xi)/\eta(t)\leq 1$ because $\eta$ decreases in $(T+r,T+2r).$
Moreover, $|\xi-T|< 3r$ for $(x',\xi)\in\Omega_r$ and
\begin{align*}
\int_0^\xi \frac{\eta(\xi)}{\eta(t)} \D_2(\xi)\tl V(x',t)dt =&\
    \int_0^T \frac{\eta(\xi)}{\eta(t)} \D_2(\xi)\tl V(x',t)dt+
     \int_T^\xi \frac{\eta(\xi)}{\eta(t)} \D_2(\xi)\tl V(x',t)dt\\
     =&\
     \eta(\xi)\int_0^T	\D_2(\xi) V(x',t)dt+
     \int_T^\xi \frac{\eta(\xi)}{\eta(t)} \D_2(\xi)\tl V(x',t)dt
\end{align*}
by means of \eqref{eqF1} and  since $\eta(t)=\eta(T)=1$ as $t\leq T.$

We get from \eqref{eq2.8f} and \eqref{eq2.9}
that $\tl V\in W^{2,p}(\Omega_r)$ solves the Dirichlet
problem
\begin{equation}\label{eq2.10}
\begin{cases}
{\L}'\tl V=F_2(x',\xi)+\ds\int_T^\xi \frac{\eta(\xi)}{\eta(t)} \D_2(\xi)\tl V(x',t)dt
	 \qquad\text{a.a.}\ (x',\xi)\in\Omega_r,\\[8pt]
\tl V=\tl\Phi:=\eta\Phi=
\begin{cases}
\eta\mu\varphi\in W^{2-1/q,q} & \text{on}\ \partial_1\Omega_r\quad
\text{(by \eqref{eqUeq})},\\
0 & \text{on}\ \partial_2\Omega_r\quad
\text{(by \eqref{eqUeq})},\\
0	   & \text{on}\ \partial\Omega_r^+\quad \text{(by \eqref{eqeta})},\\
 V\in W^{2-1/q',q'}  & \text{on}\ \partial\Omega_r^-\quad
	\text{(since $\xi<T-2r$ and}\\
&  \qquad\qquad\quad \text{$V\in W^{2,q'}(\P_{r,T})$)}
\end{cases}
\end{cases}
\end{equation}
where, recalling $V\in W^{2,q'}(\P_{r,T}),$ we have
\begin{equation}\label{eqF2}
F_2(x',\xi):= \eta F_1+\L_1 V +
\eta(\xi)    \int_0^T  \D_2(\xi) V(x',t)dt\in L^{q'}(\Omega_r).
\end{equation}
We are going to prove now that $\tl V\in W^{2,q'}(\Omega_r)$ for small enough $r>0,$ whence
it will follow $V\in W^{2,q'}(\P_{r,T+r})$ in view of
\eqref{eqeta} and $V\equiv0$ near $\partial_2\Omega_r.$
 The claim is obvious if $q'=p$ because $V\in W^{2,p}(\T_r).$
 Otherwise, take an arbitrary
$s\in[p,q']$ and denote by $W^{2,s}_*(\Omega_r)$ the Sobolev space $W^{2,s}(\Omega_r)$ normed with 
$$
\|u\|_{W^{2,s}_*(\Omega_r)}:=
\|u\|_{L^{s}(\Omega_r)}+r\|Du\|_{L^{s}(\Omega_r)}
+r^2\|D^2u\|_{L^{s}(\Omega_r)}.
$$
Define now the operator $\FF\colon W^{2,s}_*(\Omega_r)\to W^{2,s}_*(\Omega_r)$
as follows: for any $w\in W^{2,s}_*(\Omega_r)$ the image $\FF w\in
W^{2,s}_*(\Omega_r)$ is the {\it unique\/} solution of the Dirichlet problem
\begin{equation}\label{eq2.11}
\begin{cases}
{\L}'(\FF w)=F_2+\ds\int_T^\xi \frac{\eta(\xi)}{\eta(t)} \D_2(\xi) w
(x',t)dt \in L^s(\Omega_r) &\text{a.a.}\ (x',\xi)\in\Omega_r,\\[8pt]
\FF w=\tl\Phi\in W^{2-1/s,s}(\partial\Omega_r) & \text{on}\
\partial\Omega_r.
\end{cases}
\end{equation}
We will prove that  $\FF$ is a contraction for small
values of $r.$ For this goal, take arbitrary $w_1,\ w_2\in W^{2,s}_*(\Omega_r).$ The
difference $\FF w_1-\FF w_2$ solves
\begin{equation}\label{eq2.12}
\begin{cases}
{\L}'(\FF w_1-\FF w_2)=\ds
\int_T^\xi \frac{\eta(\xi)}{\eta(t)} \D_2(\xi) (w_1-w_2) (x',t)dt
 &\text{a.a.}\ (x',\xi)\in\Omega_r,\\
\FF w_1-\FF w_2=0 & \text{on}\
\partial\Omega_r.
\end{cases}
\end{equation}
In order to apply the $L^s$-{\it a~priori\/} estimates from \cite{ADN} or \cite{GT} for the
solutions of \eqref{eq2.12}, we have to control the dependence on $r$ therein. For,
we recall that $r^{-1}\Omega_r\in C^{1,1}$ uniformly in $r$ and 
apply a standard approach consisting of dilation of $\Omega_r$ onto
$r^{-1}\Omega_r,$ reduction of the problem \eqref{eq2.12} to a new one in
variables $(\tl y',\tl \xi)\in r^{-1}\Omega_r,$ application of the $L^s$-estimates from
\cite[Theorem~9.17]{GT} and finally turning back to \eqref{eq2.12}
(see the Proof of Lemma~2.2, Eq.~(2.12) in \cite{Pl2}). This way, one gets
\begin{equation}\label{eq2.13}
\|\FF w_1-\FF w_2\|_{W^{2,s}_*(\Omega_r)} \leq C
r^2 \left\|
\int_T^\xi \frac{\eta(\xi)}{\eta(t)} \D_2(\xi) (w_1-w_2)
(x',t)dt\right\|_{L^{s}(\Omega_r)}
\end{equation}
where the constant $C$ is independent of $r.$ Jensen's integral inequality yields
$$
r^2 \left\|
\int_T^\xi \frac{\eta(\xi)}{\eta(t)} \D_2(\xi) (w_1-w_2)
(x',t)dt\right\|_{L^{s}(\Omega_r)} \leq
C \max_{(x',\xi)\in\Omega_r} |\xi-T|
\|w_1-w_2\|_{W^{2,s}_*(\Omega_r)}
$$
and thus \eqref{eq2.13} rewrites into
$$
\|\FF w_1-\FF w_2\|_{W^{2,s}_*(\Omega_r)} \leq C
\max_{(x',\xi)\in\Omega_r} |\xi-T|
\|w_1-w_2\|_{W^{2,s}_*(\Omega_r)}.
$$
We have $\max_{(x',\xi)\in\Omega_r} |\xi-T|<3r,$  $C$ is independent of $r$ and therefore $\FF$ will be really a contraction from $W^{2,s}_*(\Omega_r)$ into itself
for any $s\in[p,q']$ if $r\leq r_0$ with $r_0$ under control
and small enough. Fixing $r=r_0/2,$ there is a unique fixed point of $\FF$ in $W^{2,s}_*(\Omega_r)$ for all $s\in[p,q'].$
However, $\tl V\in W^{2,p}(\Omega_r)$ is already a  fixed point of $\FF$ since it solves \eqref{eq2.10} and therefore
$\tl V\in W^{2,q'}(\Omega_r).$ It follows $V\in W^{2,q'}(\P_{r,T+r})$ by means of
$V\in W^{2,q'}(\P_{r,T}),$ $\tl V\equiv 0$ on $\T_r\sm \T_{3r/4}$ and the
properties of $\eta(\xi).$ 

\paragraph*{\it Case B: $T<T_{\max}\leq T+3r$} We have $\T_r\sm\P_{r,T}\neq\emptyset,$ $\P_{r,T+3r}\equiv\T_r$
now and we do not need anymore
the cut-off function $\eta$ because $V=\partial U/\partial\Ll\equiv 0$  near the points of $\partial_2\T_r$
where $\xi>T$ (cf. \eqref{eqMU}).  Thus,
it suffices to repeat the above arguments with $\eta(\xi)\equiv 1$ $\forall\xi\in\R$ and $\Omega_r\in C^{1,1}$ defined
as before when $\xi\leq T$ while $\T_{3r/4}\sm\P_{3r/4,T}\subset\big(\Omega_r\cap\{\xi>T\}\big)\subset \T_{r}\sm\P_{r,T}$
(cf. \eqref{eqMU}). We have anyway a problem like
\eqref{eq2.10} for $V\equiv\tl V$ with boundary condition
$$
V=\partial U/\partial\Ll=
\begin{cases}
\mu\varphi\in W^{2-1/q,q} & \text{on}\ \partial_1\Omega_r=\partial\Omega_r\cap\partial\Omega,\\
0 & \text{on}\ \partial_2\Omega_r=\partial\Omega_r\cap\Omega\cap\{\xi>T-3r\},\\
V\in W^{2-1/q',q'}  & \text{on}\ \partial\Omega_r^-\quad
	\text{(by hypothesis).}
\end{cases}
$$
Therefore, the procedure from {\it Case~A\/}  gives $V\in W^{2,q'}(\P_{r,T+3r}).$
\paragraph*{\it Case C: $T_{\max}\leq T$} We have $\P_{r,T+r}\equiv \P_{r,T}\equiv\T_r$ now and thus the claim.
\end{proof}
\begin{prp}\label{prpAEST}
Suppose $r<r_0$ with $r_0$ given in Proposition~$\ref{prpIMP}.$ Then
the solution $V$ of the problem $\eqref{eq2.8f}$ lies in $W^{2,q}(\T_{r})$ and satisfies the estimate
\begin{align}\label{eqAEST}
\|V\|_{W^{2,q}(\T_r)}\leq&\ C\Big(\|u\|_{L^q(\Omega)}
+\|f\|_{\F^q(\Omega,\N)}+\|\varphi\|_{\Phi^q(\partial\Omega,\N)}\\
\nonumber
&\qquad +\|u\|_{W^{1,q}(\T_r)}
+\|\partial u/\partial\Ll\|_{W^{1,q}(\T_r)}
\Big).
\end{align}
\end{prp}
\begin{proof}
We note that $V\in W^{2,q}\subseteq W^{2,q'}$ near $B'_r(x_0)$
in view of
$B'_r(x_0)\subset\N''\sm\N',$ Proposition~\ref{lemA1} and
\eqref{eqderDP}. Therefore, successive applications of
Proposition~\ref{prpIMP} with increasing values of $T$ will give $V\in
W^{2,q'}(\T_r),$ $q'>p.$
After that, in order to get $V\in W^{2,q}(\T_r),$
it suffices to put $q'$ in the place of $p$ in \eqref{eqq'} and
to repeat finitely many times the above arguments until
$q'=q.$    

To obtain \eqref{eqAEST}, we take $T\in(0,t^+-t^-)$ to be arbitrary, fix $r=r_0/2,$  and consider the domains
$\Omega_r$ defined in the proof of Proposition~\ref{prpIMP}. Let $\tl V=\eta V\in
W^{2,q}(\T_r)$ solve \eqref{eq2.10} with $\eta$
 given by \eqref{eqeta} in {\it Case~A\/} and $\eta\equiv 1$ in {\it Case~B.\/} Since $\tl V$
is a fixed point of the mapping $\FF\colon\ W^{2,q}(\Omega_r)\to
W^{2,q}(\Omega_r),$ $\FF \tl V=\tl V,$ we get
$$
\|D^2 \tl V\|_{L^q(\Omega_r)}=
\|D^2 (\FF \tl V) \|_{L^q(\Omega_r)}\leq
\|D^2 (\FF \tl V-\FF 0) \|_{L^q(\Omega_r)}+
\|D^2 (\FF 0) \|_{L^q(\Omega_r)},
$$
while
$$
\|D^2 (\FF w_1-\FF w_2) \|_{L^q(\Omega_r)}\leq \theta
\|D^2 (w_1-w_2) \|_{L^q(\Omega_r)}\quad \forall w_1,\ w_2\in W^{2,q}(\Omega_r),
\quad \theta<1
$$
because $\FF$ is a contraction, \eqref{eq2.13} and
the fact that $\D_2(\xi)$ is a homogeneous second-order operator (cf. \eqref{eqF1}).
This way, $\|D^2 (\FF \tl V-\FF 0) \|_{L^q(\Omega_r)}\leq \theta \|D^2 (\tl V-0) \|_{L^q(\Omega_r)}
=\theta \|D^2 \tl V \|_{L^q(\Omega_r)}$ and therefore
\begin{equation}\label{eq2.14}
\|D^2 \tl V\|_{L^q(\Omega_r)}\leq C
\|D^2 (\FF 0)\|_{L^q(\Omega_r)}
\end{equation}
with $\FF 0\in W^{2,q}(\Omega_r)$ being the unique solution of the Dirichlet problem
$$
\begin{cases}
{\L}'(\FF 0)=F_2 \quad\text{a.e.}\ \Omega_r,\qquad \FF 0=\tl\Phi \quad \text{on}\
\partial\Omega_r
\end{cases}
$$
(see \eqref{eq2.11}), for which the $L^p$-theory (cf. \cite[Chapter~9]{GT}) gives
\begin{equation}\label{eq2.15}
\|D^2 (\FF 0)\|_{L^q(\Omega_r)} \leq \|\FF 0\|_{W^{2,q}(\Omega_r)}\leq C
\left(\|F_2\|_{L^q(\Omega_r)}+
\|\tl \Phi\|_{W^{2-1/q,q}(\partial\Omega_r)}\right).
\end{equation}
Direct applications, based on \eqref{eqF2} and \eqref{eqF1}, yield
\begin{align*}
&\|F_2\|_{L^q(\Omega_r)}=
\left\|\eta F_1+\L_1V+\eta(\xi)\int_0^T \D_2(\xi) V(x',t)dt
	     \right\|_{L^q(\Omega_r)}\\
&\qquad \leq C\Big(
\|\partial F/\partial\Ll\|_{L^q(\Omega_r)}+\|U\|_{W^{2,q}(\N''\sm\N')}+\|U\|_{W^{1,q}(\T_r)}
+\|V\|_{W^{1,q}(\T_r)}\\
&\qquad\qquad\quad + \|D^2V\|_{L^{q}(\P_{r,T})}\Big)\\
&\qquad\leq C\Big(
\|\partial f/\partial\Ll\|_{L^q(\N)}+\|u\|_{W^{2,q}(\N''\sm\N')}
+\|u\|_{W^{1,q}(\T_r)}+\|\partial u/\partial\Ll\|_{W^{1,q}(\T_r)}\\
&\qquad\qquad\quad + \|D^2V\|_{L^{q}(\P_{r,T})}\Big)
\end{align*}
in view of \eqref{eq2.5}, \eqref{eqUeq}, $U=\mu u,$ $V=\partial U/\partial\Ll$ and
\eqref{eqMU}. Moreover,
\begin{align*}
\|\tl\Phi\|_{W^{2-1/q,q}(\partial\Omega_r)}
\leq&\	C\left(
\|\varphi\|_{W^{2-1/q,q}(\partial\Omega\cap\N)} +\|V\|_{W^{2,q}(\P_{r,T})}\right)\\
\leq&\	C\left(
\|\varphi\|_{W^{2-1/q,q}(\partial\Omega\cap\N)}
+\|V\|_{W^{1,q}(\T_r)}+\|D^2V\|_{L^{q}(\P_{r,T})}\right)\\
 \leq&\  C\left(
\|\varphi\|_{W^{2-1/q,q}(\partial\Omega\cap\N)}
+\|\partial u/\partial\Ll\|_{W^{1,q}(\T_r)}+\|D^2V\|_{L^{q}(\P_{r,T})}\right)
\end{align*}
by \eqref{eq2.10} and  $\partial\Omega_r^-\subset \P_{r,T}.$
Further on, $\tl V=V$ on $\P_{r,T+r},$ whence
$$
\|D^2V\|_{L^q(\P_{r,T+r})} \leq \|D^2V\|_{L^q(\P_{r,T})}+
\|D^2\tl V\|_{L^q(\Omega_{r})}.
$$
Therefore, setting $\zeta(T):=\|D^2V\|_{L^q(\P_{r,T})}$ and $K:=
\|u\|_{L^q(\Omega)}+\|f\|_{\F^{q}(\Omega,\N)}+
\|\varphi\|_{\Phi^q(\partial\Omega,\N)}+\|u\|_{W^{1,q}(\T_r)}
+\|\partial u/\partial\Ll\|_{W^{1,q}(\T_r)},$ 
it follows from \eqref{eq2.14}, \eqref{eq2.15} and Proposition~\ref{lemA1} that
\begin{equation}\label{eqzeta}
\zeta(T+r)\leq C\left(K+\zeta(T)\right)\qquad \forall T\in (0,t^+-t^-).
\end{equation}
To get \eqref{eqAEST}, we let $m$ to be the least integer such that
$T_{\max}\leq mr$ and iterate \eqref{eqzeta} in order to obtain
\begin{align*}
\|D^2V\|_{L^q(\T_{r})}=&\ \|D^2V\|_{L^q(\P_{r,T_{\max}})}=\zeta(T_{\max})=
 \zeta(mr)=\zeta((m-1)r+r)\\
\leq&\ C\big(K+\zeta((m-1)r)\big)=C\big(K+\zeta((m-2)r+r)\big)\\
\leq&\ K(C+C^2)+C^2\zeta\big((m-2)r\big)\\
\vdots &\\
\leq&\ K\sum_{j=1}^m C^j + C^m\zeta(0)= K\sum_{j=1}^m C^j
\end{align*}
This proves \eqref{eqAEST}.
\end{proof}

\begin{rem}\label{remCONST} \em
It is important to note that the constant $C$ in Proposition~$\ref{prpAEST}$ depends on
$m$ through $T_{\max},$ and therefore on the point $x_0\in \E.$ Actually, that constant
will have the very same value for each other point of $\E$ lying on the same $\Ll$-trajectory
as $x_0.$

Moreover, if the {\it improving-of-integrability property\/} asserted in 
Propositions~$\ref{prpIMP}$ and $\ref{prpAEST}$ holds on a set $S\subset\ol\Omega$ then it is
guaranteed, on the base of \eqref{eqIFF}, on any other set which can be reached from $S$ along
$\Ll$-trajectories.
\end{rem}
To complete the proof of Lemma~\ref{lemA2}, we 
select a finite set $\{\T^j_r\}_{j=1}^N$ of neighbourhoods covering the {\it compact\/} $\E,$ each of the type $\T_r$ above with $r=r_0/2,$
and such that $\T:=\mathrm{closure\,}\left(\bigcup_{j=1}^N \T^j_{r/2}\right)\subset\N''$ is a closed neighbourhood of $\E$
in $\ol\Omega.$ It is clear that Proposition~\ref{lemA1} remains true with $\T$ instead of $\N'$ and then
\eqref{eq2.5} rewrites into
\begin{equation}\label{eq2.5'}
\|u\|_{W^{2,q}(\Omega\sm\T)} \leq C\left(
\|u\|_{L^{q}(\Omega)}+
\|f\|_{L^{q}(\Omega)}+
\|\varphi\|_{W^{1-1/q,q}(\partial\Omega)}\right).
\end{equation}
The
{\it improving-of-integrability\/} claimed in Lemma~\ref{lemA2} then follows from
\eqref{eqIFF}, Proposition~\ref{prpAEST} and \eqref{eq2.5'} (recall $U=u$ on $\T^j_{r/2}$).
Similarly, \eqref{eqIFF}, \eqref{eq2.5'} and \eqref{eqAEST} yield
\begin{align}\label{eqAPPR}
\|u\|_{W^{2,q}(\N'')}\leq&\ \|u\|_{W^{2,q}(\T)}+\|u\|_{W^{2,q}(\N''\sm\T)}\\
\nonumber \leq&\ 
C\big( \|u\|_{L^{q}(\Omega)}
+\|f\|_{\F^{q}(\Omega,\N)}+ \|\varphi\|_{\Phi^q(\partial\Omega,\N)}\\
\nonumber
&\qquad +\|u\|_{W^{1,q}(\N)} +\|\partial u/\partial\Ll\|_{W^{1,q}(\N)} \big).
\end{align}
Later on, $\N\sm\N''\subset \Omega\sm\N'$ and
\begin{align*}
\|u\|_{W^{1,q}(\N)}\leq&\ \|u\|_{W^{1,q}(\N'')} +\|u\|_{W^{1,q}(\N\sm\N'')}\\
\leq&\ \varepsilon \|u\|_{W^{2,q}(\N'')}+C(\varepsilon)\big(\|u\|_{L^{q}(\Omega)}+\|u\|_{W^{2,q}(\Omega\sm\N')}\big)
\end{align*}
in view of the interpolation inequality for the $W^{2,q}(\N'')$-norms with $\varepsilon>0$ under control\footnote{This requires 
some minimal smoothness of $\partial\N''$ and it is not restrictive to take it Lipschitz continuous at the very beginning.}.
In the same manner,
\begin{align*}
\|\partial u/\partial\Ll\|_{W^{1,q}(\N)}\leq&\ \|\partial u/\partial\Ll\|_{W^{1,q}(\N')}
  +\|\partial u/\partial\Ll\|_{W^{1,q}(\N\sm\N')}\\
\leq&\ \varepsilon \|\partial u/\partial\Ll\|_{W^{2,q}(\N')}+
  C(\varepsilon)\big(\|\partial u/\partial\Ll\|_{L^{q}(\N')}+\|u\|_{W^{2,q}(\Omega\sm\N')}\big),
\end{align*}
while
$$
\|\partial u/\partial\Ll\|_{W^{2,q}(\N')} \leq C\big(\|u\|_{W^{2,q}(\N'')}+
\|u\|_{L^{q}(\Omega)} + \|f\|_{\F^{q}(\Omega,\N)}+ \|\varphi\|_{\Phi^q(\partial\Omega,\N)}
\big)
$$
by means of the local {\it a~priori\/} estimates (\cite[Theorem~9.11]{GT}) for the problem \eqref{eqderDP}.

A substitution of the above expressions into \eqref{eqAPPR} and \eqref{eq2.5} give
\begin{align*}
\|u\|_{W^{2,q}(\N'')}\leq&\ C\big( \|u\|_{L^{q}(\Omega)}
+\|f\|_{\F^{q}(\Omega,\N)}+ \|\varphi\|_{\Phi^q(\partial\Omega,\N)}\\
&\qquad +\varepsilon\|u\|_{W^{2,q}(\N'')} +C(\varepsilon)\|\partial u/\partial\Ll\|_{L^{q}(\N')} \big)
\end{align*}
whence, choosing $\varepsilon>0$ small enough, we get
$$
\|u\|_{W^{2,q}(\N'')}\leq C\big( \|u\|_{L^{q}(\Omega)}
+\|f\|_{\F^{q}(\Omega,\N)}+ \|\varphi\|_{\Phi^q(\partial\Omega,\N)}
+\|u\|_{W^{1,q}(\N')} \big).
$$
Similarly, another application of the interpolation inequality yields
$$
\|u\|_{W^{1,q}(\N')}\leq \|u\|_{W^{1,q}(\N'')}\leq \delta \|u\|_{W^{2,q}(\N'')}+C(\delta)\|u\|_{L^{q}(\N'')}
$$
and thus 
$$
\|u\|_{W^{2,q}(\N'')}\leq C\big( \|u\|_{L^{q}(\Omega)}
+\|f\|_{\F^{q}(\Omega,\N)}+ \|\varphi\|_{\Phi^q(\partial\Omega,\N)}\big).
$$
for small $\delta>0.$
 The proof of Lemma~\ref{lemA2} is completed.
\end{proof}
The statement of Theorem~\ref{thmINC} follows from
Proposition~\ref{lemA1} and Lemma~\ref{lemA2}.

\end{document}